\documentclass[preprint,12pt,authoryear]{elsarticle}

\usepackage[left=2.2cm,right=2.2cm,top=2.5cm,bottom=2.5cm]{geometry}
\usepackage{lmodern}
\usepackage[T1]{fontenc}
\usepackage[utf8]{inputenc}
\usepackage{amsmath,amssymb,amsthm}
\usepackage{mathtools}
\usepackage{bm}
\usepackage{graphicx}
\usepackage{subcaption}
\usepackage{booktabs}
\usepackage{multirow}
\usepackage{array}
\usepackage{enumitem}
\usepackage{hyperref}
\usepackage[capitalise,noabbrev]{cleveref}
\usepackage{tikz}
\usetikzlibrary{arrows.meta,positioning,shapes.geometric,fit,backgrounds,calc}
\usepackage{pgfplots}
\pgfplotsset{compat=1.18}
\usepackage{algorithm}
\usepackage{algpseudocode}

\theoremstyle{definition}

\journal{Transportation Research Part C (preprint)}

\begin{document}

\begin{frontmatter}

\title{Route Based Map Matching \\ via a Structured Codebook and Token Sequence Decoding}

\author[aff1]{Takara Sakai\corref{cor1}}
\ead{sakai.t.dcad@m.isct.ac.jp}
\cortext[cor1]{Corresponding author.}
\affiliation[aff1]{organization={Department of Civil and Environmental Engineering,\\ Institute of Science Tokyo},
                   addressline={2-12-1 W6-9, Ookayama, Meguro, Tokyo},
                   postcode={152-8550},
                   country={Japan}}

\author[aff2]{Hiroyuki Hasada}
\affiliation[aff2]{organization={Institute of Industrial Science,\\ The University of Tokyo},
                   addressline={Komaba 4-6-1, Meguro-ku, Tokyo},
                   postcode={153-8505},
                   country={Japan}}

\begin{abstract}
  This study proposes an efficient and computationally light route based map matching method for GPS track data on urban expressway networks.
  The key idea is to exploit a symbolic structure of named lines and named junctions that link level map matching leaves unused.
  We represent each candidate route as a sequence of line and junction names, take the set of such sequences as a route codebook, and formulate map matching as scored alignment of a probe trajectory against members of the codebook.
  Probes become token sequences via a mesh quantizer, a precomputed grid mapping each coordinate to a line or junction token, and the decoder returns a member of the codebook by construction.
  The codebook is indexed by a DAFSA $\times$ Levenshtein automaton, a fuzzy lookup technique from approximate string matching and speech recognition; the per query decoding cost is orders of magnitude lower than a brute force scan.
  We evaluate the method on a deformed replica of the Tokyo Metropolitan Expressway topology.
  The method recovers the exact route at moderate GPS noise and continues to identify the line and junction sequence under heavy noise; a sensitivity analysis maps the mesh resolution operating range.
  Real probe evaluation, channel model calibration, and a head to head HMM comparison are left to a forthcoming version.
\end{abstract}


\end{frontmatter}

\section{Introduction}
\label{sec:intro}

Map matching, the problem of associating a noisy GPS trajectory with a path on a road network, has been dominated for over a decade by per link, per observation decoding: each observation is matched to the most likely link, and the link sequence is recovered by Viterbi decoding \citep{Newson2009}.
The broader landscape is surveyed by \citet{Quddus2007}, and speed oriented descendants such as fmm \citep{Yang2018} amortize the Viterbi cost while keeping the per link decoding granularity.
\par
However, on urban expressway networks per link decoding faces two structural difficulties.
First, a single misidentification at a junction propagates into a wrong route for the remainder of the trajectory, since the choice of branch determines all downstream links.
Second, main lines, branching ramps, and elevated decks of a ring road and a radial route sit within a few metres of each other on a multilayer structure, so GPS noise easily renders alternative links indistinguishable.
\par
To address both difficulties, this paper shifts the unit of map matching from links to routes.
The decoding decision is then made globally over a precomputed set of valid routes, rather than locally at each observation, which prevents an early junction error from propagating downstream.
By construction the decoder returns a route that exists in the candidate set; a per link decoder makes independent per observation decisions and can produce a link sequence that corresponds to no valid route.
This reformulation rests on a symbolic structure that general road networks lack: every segment of an urban expressway belongs to a named line, and every junction and interchange has a proper name.
A driver describes a trip as ``took C1 inbound, transferred at Hamazakibashi JCT to Route~3, exited at Shibuya'', not as a list of link IDs.
\par
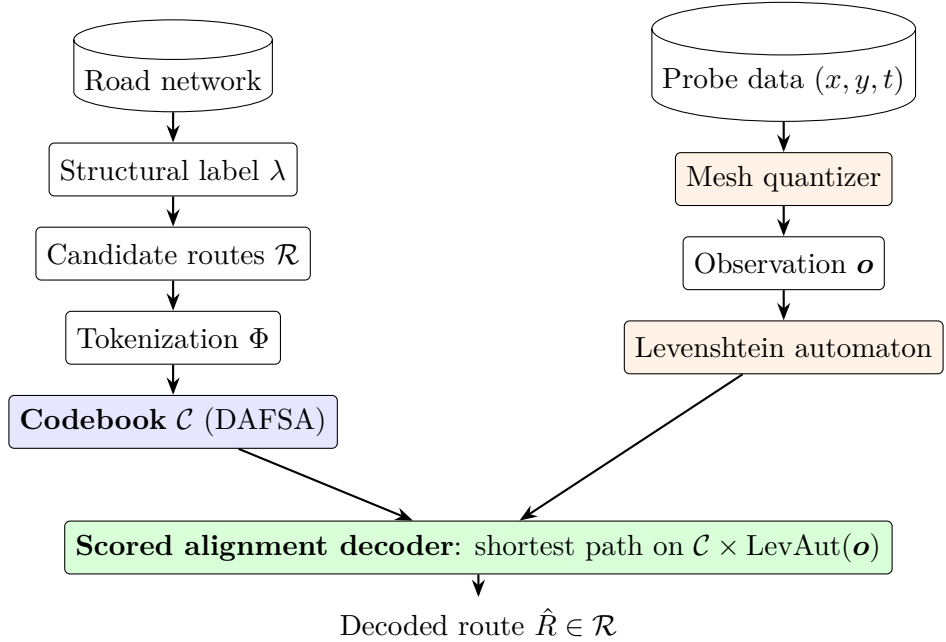
\begin{figure}[tbp]
  \centering
  \begin{tikzpicture}[
    >=Stealth,
    every node/.style={font=\small},
    box/.style={rectangle, draw, rounded corners=2pt, align=center, minimum height=7mm, inner sep=4pt},
    data/.style={cylinder, draw, shape border rotate=90, aspect=0.2, minimum height=6mm, minimum width=20mm, align=center},
    arrow/.style={->, thick},
  ]
    \node[data] (net) {Road network};
    \node[box, below=4mm of net] (label) {Structural label $\lambda$};
    \node[box, below=4mm of label] (routes) {Candidate routes $\mathcal{R}$};
    \node[box, below=4mm of routes] (phi) {Tokenization $\Phi$};
    \node[box, below=4mm of phi, fill=blue!10] (codebook) {\textbf{Codebook} $\mathcal{C}$ (DAFSA)};

    \node[data, right=50mm of net] (probe) {Probe data $(x,y,t)$};
    \node[box, below=4mm of probe, fill=orange!10] (mesh) {Mesh quantizer};
    \node[box, below=4mm of mesh] (obs) {Observation $\bm{o}$};
    \node[box, below=4mm of obs, fill=orange!10] (lev) {Levenshtein automaton};

    \node[box, below=18mm of $(codebook)!0.5!(lev)$, fill=green!15, minimum width=75mm] (dec) {\textbf{Scored alignment decoder}: shortest path on $\mathcal{C} \times \mathrm{LevAut}(\bm{o})$};
    \node[below=3mm of dec] (out) {Decoded route $\hat{R} \in \mathcal{R}$};

    \draw[arrow] (net) -- (label);
    \draw[arrow] (label) -- (routes);
    \draw[arrow] (routes) -- (phi);
    \draw[arrow] (phi) -- (codebook);
    \draw[arrow] (probe) -- (mesh);
    \draw[arrow] (mesh) -- (obs);
    \draw[arrow] (obs) -- (lev);
    \draw[arrow] (codebook) -- (dec);
    \draw[arrow] (lev) -- (dec);
    \draw[arrow] (dec) -- (out);
  \end{tikzpicture}
  \caption{Processing flow. The route side becomes a codebook of token sequences; the probe side becomes an observation sequence via mesh lookup; the decoder returns the best aligned codebook member.}
  \label{fig:framework}
\end{figure}

A finite vocabulary $\mathcal{T}$ collects line, junction, and interchange tokens from each link, and a candidate route maps to the token sequence its links carry, with consecutive duplicates contracted.
The set of these sequences is the route \emph{codebook} $\mathcal{C}$.
A probe trajectory becomes an observation token sequence via a mesh quantizer: the network is rasterized once into a grid in which each cell carries the line or junction token it covers, and each probe point is mapped to its token by integer division and table lookup, with no per query link distance computation at inference time.
Decoding is the scored alignment of the observation against members of $\mathcal{C}$.
\cref{fig:framework} shows the processing flow.

For the search to scale, we index the codebook with finite state automata.
A prefix sharing trie \citep{Fredkin1960} and the deterministic acyclic finite state automaton (DAFSA) of \citet{Daciuk2000} represent the codebook compactly through prefix and suffix sharing, and the product of the DAFSA with the weighted Levenshtein automaton of \citet{Schulz2002} casts bounded edit fuzzy lookup as a shortest path search on the product automaton.
The same product construction is a standard tool in approximate string matching and is the dictionary backbone of weighted finite state transducer (WFST) decoders in speech recognition \citep{Mohri2002}.
On the route codebook the index lowers the per query decoding cost by orders of magnitude relative to a brute force scan, while symbolic trajectory work such as SAX \citep{Lin2007} operates at different granularities (time series) and does not target network level route identity.

The contributions of this paper are the following two.
\begin{enumerate}[leftmargin=1.5em]
  \item We cast urban expressway map matching as scored alignment over a route codebook whose vocabulary consists of line and junction names, with a mesh quantizer that maps each probe coordinate to a line or junction token by table lookup, replacing per point link distance computation.
  \item We compare prefix sharing and minimal automaton indices for the codebook and combine the latter with a Levenshtein automaton for fuzzy decoding; the combined index lowers the per query decoding cost by orders of magnitude relative to a brute force scan.
\end{enumerate}

This version focuses on the methodological formulation and a controlled synthetic validation.
The evaluation uses a toy network data set: a topological replica of the Tokyo Metropolitan Expressway (760 links, 695 nodes, 15 named lines, 24 junctions, 102 interchanges) with deliberately deformed node coordinates, and synthetic probe trajectories generated from ground truth routes with an additive Gaussian noise model.
Validation on real probe data, channel model calibration, and a head to head comparison with an HMM baseline are reserved for the journal version.
The remainder of the paper is organized as follows: \cref{sec:method} presents the methodology, \cref{sec:setup} describes the experimental setup, \cref{sec:results} reports results across three experiments (noise robustness, mesh resolution sensitivity, and indexing structure efficiency), and \cref{sec:conclusion} concludes.

\section{Methodology}
\label{sec:method}

The methodology follows the processing flow of \cref{fig:framework}.
The route side (\cref{sec:method_language}) turns candidate routes into token sequences over a line name alphabet and forms the codebook $\mathcal{C}$.
The probe side (\cref{sec:method_mesh} and \cref{sec:method_channel}) turns each probe trajectory into an observation token sequence via a mesh quantizer, without per query link distance computation at inference time.
The decoder (\cref{sec:method_coding}) returns the member of $\mathcal{C}$ with the highest alignment score.
\cref{sec:method_distance} discusses pairwise separation in $\mathcal{C}$ and the conditions under which $\Phi$-collisions can occur, and \cref{sec:method_index} develops the index structures on which decoding is executed.

\begin{figure}[tbp]
  \centering
  \begin{subfigure}[t]{0.32\textwidth}
    \centering
    \includegraphics[width=\textwidth]{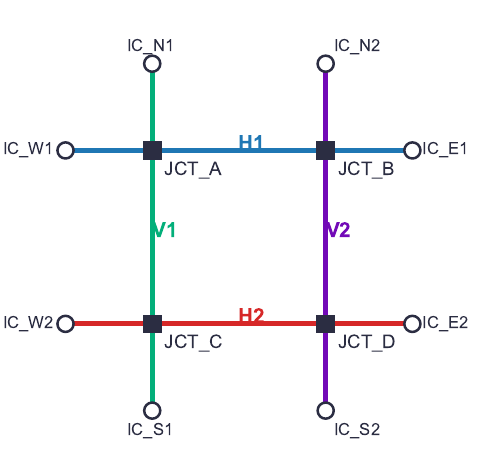}
    \caption{Network.}
    \label{fig:igrid_network}
  \end{subfigure}
  \hfill
  \begin{subfigure}[t]{0.32\textwidth}
    \centering
    \includegraphics[width=\textwidth]{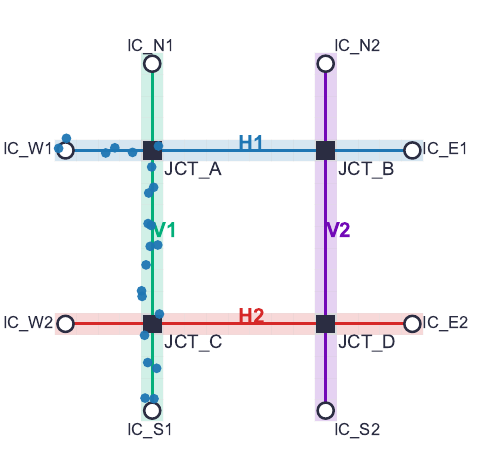}
    \caption{Mesh quantizer with noisy probes.}
    \label{fig:igrid_mesh}
  \end{subfigure}
  \hfill
  \begin{subfigure}[t]{0.32\textwidth}
    \centering
    \includegraphics[width=\textwidth]{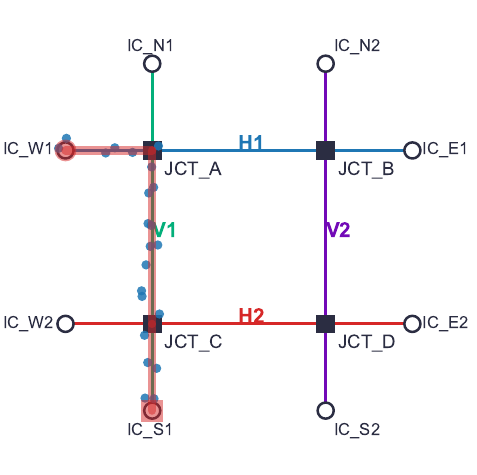}
    \caption{Decoded route on top of the probes.}
    \label{fig:igrid_match}
  \end{subfigure}
  \caption{The \#-shaped running example. Two horizontal lines (H1, H2) and two vertical lines (V1, V2) cross at four named junctions \texttt{JCT\_A--D}, with eight named interchanges at the tips (squares: JCT, circles: IC). Panel (b) shows the mesh quantizer applied to a noisy probe sequence along the route \texttt{IC\_W1 $\to$ JCT\_A $\to$ JCT\_C $\to$ IC\_S1}; panel (c) shows the decoded route, recovered as the codeword \texttt{[IN@W1, LINE\_H1, JCT\_A, LINE\_V1, OUT@S1]}.}
  \label{fig:igrid}
\end{figure}

We illustrate the three stages on a small synthetic \#-shaped network running throughout this section.
The network has two horizontal lines (H1, H2) and two vertical lines (V1, V2) crossing at four named junctions \texttt{JCT\_A--D}, with eight named interchanges at the line tips (\cref{fig:igrid_network}).
On this network we trace a single route from \texttt{IC\_W1} through \texttt{JCT\_A} and \texttt{JCT\_C} to \texttt{IC\_S1}: \cref{fig:igrid_mesh} shows the mesh quantizer applied to a noisy probe sequence sampled along this route, and \cref{fig:igrid_match} shows the decoded route overlaid on the same probes.

\subsection{Routes as codewords: a line name language}
\label{sec:method_language}

Each link of an urban expressway network carries structural semantic labels: a line name, a junction name, and the name of the interchange where it enters or exits the expressway.
From these labels we construct a finite vocabulary $\mathcal{T}$ organized in four categories: entry tokens \texttt{IN@*}, exit tokens \texttt{OUT@*}, line tokens \texttt{LINE\_*}, and junction tokens \texttt{JCT\_*}.
A candidate route $R = (\ell_{1}, \dots, \ell_{m})$ is mapped to a token sequence by the normalization
\begin{align}
  \Phi(R) = \mathrm{canonicalize}\bigl(\lambda(\ell_{1}), \dots, \lambda(\ell_{m})\bigr) \in \mathcal{T}^{*},
  \label{eq:phi}
\end{align}
which concatenates the structural label of each link and contracts consecutive duplicates.
On the \#-shaped running example, the route \texttt{IC\_W1 $\to$ JCT\_A $\to$ JCT\_C $\to$ IC\_S1} of \cref{fig:igrid_match} produces the codeword \texttt{[IN@W1, LINE\_H1, JCT\_A, LINE\_V1, OUT@S1]}.
\par
The labelling carries junction identity at the node level rather than the link level: a node whose name ends with \texttt{JCT} marks the junction it belongs to, and $\lambda$ emits a \texttt{JCT\_<name>} token on any link with such a node as an endpoint.
Two physically distinct traversals of the same ring (for example clockwise versus counter clockwise on C1, or two passes of the C1 mainline that visit different junctions) then yield distinct codewords through their interior junction sequences, even though the line token alone would be the same.
Decoding operates entirely in this token space; the decoder consumes no link level geometric information.

\subsection{Scored alignment over the route codebook}
\label{sec:method_coding}

The image of the candidate route set $\mathcal{R}$ under the normalization $\Phi: \mathcal{R} \to \mathcal{T}^{*}$ is the route codebook
\begin{align}
  \mathcal{C} := \{\Phi(R) \mid R \in \mathcal{R}\} \subset \mathcal{T}^{*}.
  \label{eq:codebook}
\end{align}
We refer to elements of $\mathcal{C}$ as ``codewords'' and to $\mathcal{T}$ as the ``alphabet''; the vocabulary is borrowed for organizational convenience and not a claim that $\mathcal{C}$ is a designed error correcting code.
The map $\Phi$ is in general not injective: distinct routes can collapse to the same token sequence, and the decoder then returns an equivalence class of routes rather than a single one.
The frequency of such collisions is an empirical property of the network and the labelling; it is zero by construction on the toy network used here, and \cref{sec:method_distance} examines the conditions under which it can be non-zero.

Let $\bm{o} \in \mathcal{T}^{*}$ denote the observation token sequence obtained from a probe trajectory by the mesh quantizer of \cref{sec:method_mesh}.
For each $\bm{c} \in \mathcal{C}$, we compute an alignment score $S(\bm{o}, \bm{c})$ by dynamic programming with match, observation skip, and route skip operations (\cref{sec:method_channel}).
The decoded route is
\begin{align}
  \hat{R}(\bm{o}) := \Phi^{-1}\!\left(\arg\max_{\bm{c} \in \mathcal{C}} S(\bm{o}, \bm{c})\right),
  \label{eq:decoder}
\end{align}
which, by construction, is a member of $\mathcal{R}$.
A per link decoder that makes independent decisions at each observation time does not have this property; its output can correspond to no member of $\mathcal{R}$.
A direct empirical comparison with such a decoder is left to the journal version.

\subsection{Token level noisy channel model}
\label{sec:method_channel}

The proposed decoder operates entirely in token space and performs no per query link distance computation at inference time.
A probe trajectory is treated as a noisy received word, and decoding is performed by dynamic programming alignment under the token confusion matrix.

For each point of a probe trajectory, the mesh cell containing the point determines a token (a line, junction, or interchange label), and the per cell sequence is contracted by run length encoding on the dominant token to produce an observation token sequence $\bm{o} = (o_{1}, \dots, o_{n})$.
The mesh is rasterized once offline from the network's link geometry; at inference time the line assignment consists of the integer division that locates the cell and a single table lookup, with no per query link distance computation.
The observation sequence contains all token categories that the mesh can attest spatially: \texttt{IN@*}, \texttt{LINE\_*}, \texttt{JCT\_*}, and \texttt{OUT@*}.

The confusion probability $P(o \mid t)$, the probability of observing token $o$ when the true token is $t$, is defined from the topological adjacency of lines:
\begin{align}
  P(o \mid t) \propto \exp\!\left(-\frac{d_{\mathrm{topo}}(o, t)^{2}}{2 \sigma_{\mathrm{conf}}^{2}}\right),
  \label{eq:confusion}
\end{align}
where $d_{\mathrm{topo}}(o, t)$ is the hop distance between the lines of tokens $o$ and $t$ on the line adjacency graph, and $\sigma_{\mathrm{conf}}$ is a sharpness parameter.
The confusion matrix depends only on graph adjacency, not on coordinate information.
Tokens on the same line have distance zero; tokens on adjacent lines have small distance; distant lines have large distance.
We caution that this $P(o \mid t)$ is not a calibrated physical observation model derived from a noise process; it should be read as a token space smoothing kernel that softens the alignment score around adjacent lines and is a parameter of the decoder, set by hand here and learnable from data in future work.

The alignment score of a route token sequence $\Phi(R) = (r_{1}, \dots, r_{m})$ and an observation token sequence $\bm{o} = (o_{1}, \dots, o_{n})$ is computed by dynamic programming with three operations.
\begin{itemize}[leftmargin=1.5em]
  \item Match: when observation $o_{i}$ is matched to route token $r_{j}$, the score is incremented by $\log P(o_{i} \mid r_{j})$.
  \item Observation skip (insertion noise): when observation $o_{i}$ is not matched to any route token, a penalty proportional to its information cost is charged.
  \item Route skip: when route token $r_{j}$ is not matched to any observation, a penalty proportional to its information cost is charged.
\end{itemize}
The decoded route is obtained by maximizing the alignment score over all candidates, as in \cref{eq:decoder}.

A naive implementation evaluates every candidate route and runs in time $O(|\mathcal{R}| \cdot L_{R} \cdot L_{\mathrm{obs}})$.
In the proposed implementation, the token sequences of all candidates are indexed in a trie, a shared prefix tree, and the DP alignment is executed once on the trie.
Routes with a common prefix share trie nodes, so the DP is evaluated without duplication.
The total time is $O(|\mathrm{trie}| \cdot L_{\mathrm{obs}})$, where $|\mathrm{trie}|$ is the number of trie nodes; prefix sharing makes this number substantially smaller than $|\mathcal{R}| \cdot L_{R}$.

\subsection{Pairwise separation in the codebook}
\label{sec:method_distance}

The difficulty of discrimination inside $\mathcal{C}$ depends on how close its members are as token sequences.
We measure pairwise distance in $\mathcal{C}$ by the Hamming distance between equal length codewords,
\begin{align}
  d_{H}(\bm{c}_{1}, \bm{c}_{2}) := \lvert \{ i : (c_{1})_{i} \neq (c_{2})_{i} \} \rvert,
  \label{eq:dmin}
\end{align}
and the edit distance for unequal length codewords.
The minimum $d_{H}$ over distinct $\bm{c}_{1}, \bm{c}_{2} \in \mathcal{C}$ we denote $d_{\min}(\mathcal{C})$.
On the test networks of this paper $d_{\min}(\mathcal{C}) = 1$: some route pairs differ only at a single transfer junction token (concrete example in \cref{sec:res_coding}).
The hard decision argument that a code with $d_{\min} = 1$ admits no guaranteed error correction applies here, but the decoder of \cref{eq:decoder} is not hard decision: it accumulates the scored alignment of \cref{sec:method_channel} over full sequences, and continuous differences in score resolve pairs that the hard decision argument would not distinguish.
The experimental question addressed in \cref{sec:res_acc} is how well this continuous scoring discriminates on specific networks and noise levels.
\par
The $\Phi$-injectivity assumption of \cref{sec:method_coding} does not hold on every network.
When $\Phi(R_{1}) = \Phi(R_{2})$ for distinct routes $R_{1}, R_{2}$, the decoder returns the equivalence class of routes that share the same token sequence, not a single route.
For downstream tasks that are invariant under this equivalence (travel time aggregation on shared segments, for example), the ambiguity is harmless; for tasks that require the underlying link sequence, the ambiguity must be resolved by auxiliary information not modeled here.
The frequency of $\Phi$-collisions on a given network is an empirical quantity; on the networks used in \cref{sec:setup_nets} the collision rate is zero by construction.

\subsection{Mesh quantizer}
\label{sec:method_mesh}

The conversion from probe coordinates to tokens can be read as a scalar quantizer in the sense of signal processing and communication theory.
Link based map matching quantizes a probe point to a link; the proposed method quantizes into the line space instead.

The road network is rasterized at a resolution $\delta$, producing a two dimensional grid $\mathcal{M} = \{(i, j)\}$.
For each cell we accumulate, by densely sampling every link's polyline, the multiset of tokens that pass through it (one of \texttt{LINE\_*}, \texttt{JCT\_*}, \texttt{IN@*}, \texttt{OUT@*}, depending on the link's labelling).
The cell label
\begin{align}
  L(i, j) \;=\; \arg\max_{t} \, n_{ij}(t),
  \label{eq:celllabel}
\end{align}
is the most frequent token in the cell; ties are broken in the priority order \texttt{JCT\_*} $>$ \texttt{IN@*} $>$ \texttt{OUT@*} $>$ \texttt{LINE\_*}, with the lexicographically smallest token chosen within a category, and cells that no link visits are assigned the empty token $\emptyset$.
The quantization of a probe point $(x, y)$ is
\begin{align}
  \mathrm{Quant}(x, y) = L(\lfloor x / \delta \rfloor, \lfloor y / \delta \rfloor),
  \label{eq:quantizer}
\end{align}
and probes that hit empty cells are dropped from the observation sequence (their position is preserved by the surrounding probes).
The vocabulary $\mathcal{T}_{\mathrm{cell}} = \{\texttt{LINE\_*}, \texttt{JCT\_*}, \texttt{IN@*}, \texttt{OUT@*}, \emptyset\}$ is the cell label space; consecutive cells with identical labels are run length contracted on the way to the observation sequence.
\cref{fig:igrid_mesh} illustrates the mesh structure on the \#-shaped running example: each cell is shaded by its dominant token, and a probe point is mapped to a token by integer division and table lookup with no per point iteration over candidate links.

Link based map matching searches over all links at each probe point,
\begin{align}
  \hat{\ell} = \arg\min_{\ell \in \mathcal{L}} \lVert (x, y) - m(\ell) \rVert, \notag
\end{align}
at $O(|\mathcal{L}|)$ cost per point.
The proposed mesh quantization $\hat{t} = L(\lfloor x / \delta \rfloor, \lfloor y / \delta \rfloor)$ runs in $O(1)$.

The mesh resolution $\delta$ governs a tradeoff between quantization error at cell boundaries and robustness to GPS noise of scale $\sigma$.
A larger $\delta$ yields more stable line assignment near boundaries but increases cell level confusion.
A smaller $\delta$ reduces boundary sensitivity but increases the length of the observation sequence.
When $\delta \gg \sigma$, most probe points fall in cells corresponding to the true line, and confusion is low.
When $\delta \approx \sigma$, confusion increases near cell boundaries.
When $\delta \ll \sigma$, confusion approaches uniform and the quantizer carries negligible information.
Experiment 2 (\cref{sec:res_mesh}) traces this tradeoff empirically across three categorical resolution settings.

\subsection{Codebook indexing}
\label{sec:method_index}

The codebook is indexed by three sequence level structures that share the same scored alignment rule but trade space against time.
The trie of \cref{sec:method_channel} shares prefixes; the DAFSA obtained by backward suffix merging \citep{Daciuk2000} adds suffix sharing and is the minimal deterministic automaton recognizing $\mathcal{C}$; the product of the DAFSA with the weighted Levenshtein automaton of an observation $\bm{o}$ \citep{Schulz2002} casts decoding as a shortest path.
\par
The Levenshtein automaton of $\bm{o}$ with error bound $k$ accepts exactly those sequences within edit distance $k$.
Its weighted variant carries the costs of \cref{sec:method_channel}: match $-\log P(o \mid t)$, substitution from the confusion matrix, and observation/route skips at the per token cost set in \cref{sec:method_channel}.
The product automaton $A = \mathrm{DAFSA} \times \mathrm{LevAut}(\bm{o})$ has states $(q_{D}, q_{L})$, and a source to sink shortest path by Dijkstra returns the codeword of smallest total cost in $O(|A| \log |A|)$.
Prefix and suffix sharing in the DAFSA combined with the bounded fanout of the Levenshtein automaton keeps $|A|$ substantially smaller than $|\mathcal{C}| \cdot L_{\mathrm{obs}}$; the same product construction is the dictionary backbone of WFST decoders in speech recognition \citep{Mohri2002}.
The DAFSA $\times$ Levenshtein search is therefore a budgeted decoder: codewords whose alignment cost exceeds the edit budget $k$ are pruned at the product level rather than scored.
At a moderate budget the optimum can be missed when the true codeword lies just outside the budget; in this paper we set $k$ to a value at which Experiment 1 shows the proposed decoder remaining at the trie's accuracy, and we report the speed up of the budgeted decoder against the same scoring rule run by an unbudgeted scan or trie traversal.
\cref{sec:res_index} benchmarks the three sequence level indices on toy network data.

\section{Experimental Setup}
\label{sec:setup}

\subsection{Test network: toy network data}
\label{sec:setup_nets}

The evaluation testbed is the toy network data introduced in \cref{sec:intro}: the Tokyo Metropolitan Expressway (Shuto) topology with deliberately deformed node coordinates.
The deformation preserves the link, line, and junction structure but alters local distances and angles, both for visibility and to make the symbolic content of the framework visible without distraction from real geographic detail.
The graph has 695 nodes, 760 directed links, two ring routes (C1, C2) and 13 further routes (\cref{tab:lines}), 24 junctions (JCTs), and 102 interchanges (ICs), with 118 entry ramps and 118 exit ramps.
Junction identity is carried at the node level: a node whose name ends with \texttt{JCT} marks the junction it belongs to, and the labelling map $\ell$ emits a \texttt{JCT\_<name>} token on any link with such a node as an endpoint.
This treatment lets a route through C1 mainline that passes a junction without taking the branch still produce the corresponding JCT token, which is necessary to distinguish inner loop and outer loop traversals of the C1 ring through their junction sequences.
\par
\begin{figure}[tbp]
  \centering
  \includegraphics[width=0.92\textwidth]{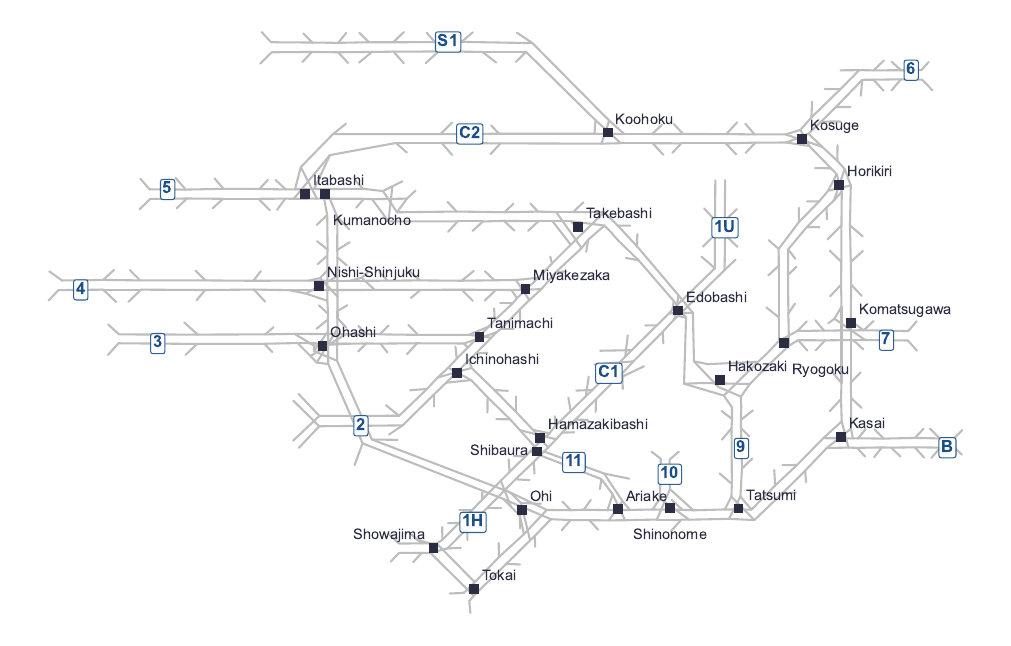}
  \caption{The toy network data: Tokyo Shuto Expressway topology with deliberately deformed coordinates. Squares: 24 junctions with romanised labels. Boxed labels (C1, C2, B, 1H, 3, $\dots$) mark each of the 15 named lines listed in \cref{tab:lines}.}
  \label{fig:network}
\end{figure}
\begin{figure}[tbp]
  \centering
  \includegraphics[width=0.92\textwidth]{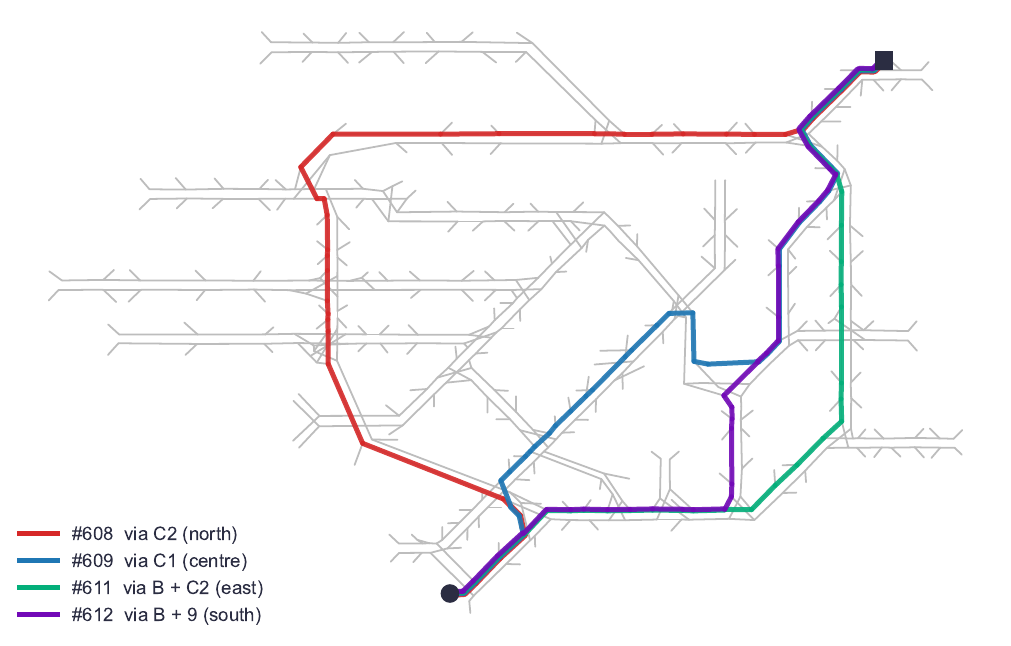}
  \caption{Within-OD route variety on \texttt{IN@Kuko-chuo $\to$ OUT@Yashio}: four representative codewords (north via C2, centre via C1, east via B+C2, south via B+9) sharing many endpoints but differing in interior line and junction sequence.}
  \label{fig:route_variety}
\end{figure}

Candidate routes are enumerated for every reachable origin destination pair (entry IC, exit IC) in the directed graph using a \emph{shortest path plus slack} cutoff: for each pair $(s, t)$ with shortest path length $\ell^{\star}(s, t)$ in links, all simple paths with at most $\ell^{\star}(s, t) + 10$ links are enumerated, and physical paths sharing the same line name codeword $\Phi(r)$ are deduplicated by keeping the shortest length representative.
\cref{tab:network_stats} summarizes the resulting codebook statistics; \cref{fig:network} visualizes the topology, and \cref{fig:route_variety} illustrates the typical within-OD route variety on the OD pair \texttt{IN@Kuko-chuo $\to$ OUT@Yashio} (Kuko-chuo IC to Yashio IC).

\begin{table}[tbp]
  \caption{Named lines of the toy network data and their formal Tokyo Shuto Expressway names.}
  \label{tab:lines}
  \centering
  \small
  \begin{tabular}{l l l}
\toprule
Token & Short name & Full name \\
\midrule
\multicolumn{3}{l}{\textit{Ring routes}} \\
\midrule
\texttt{C1} & C1 & Inner Circular Route \\
\texttt{C2} & C2 & Central Circular Route \\
\midrule
\multicolumn{3}{l}{\textit{Routes}} \\
\midrule
\texttt{B}  & B  & Wangan Route \\
\texttt{S1} & S1 & Saitama Route 1 (Kawaguchi) \\
\texttt{1H} & 1H & Route 1 Haneda Line \\
\texttt{1U} & 1U & Route 1 Ueno Line \\
\texttt{2}  & 2  & Route 2 Meguro Line \\
\texttt{3}  & 3  & Route 3 Shibuya Line \\
\texttt{4}  & 4  & Route 4 Shinjuku Line \\
\texttt{5}  & 5  & Route 5 Ikebukuro Line \\
\texttt{6}  & 6  & Route 6 Mukojima / Misato Line \\
\texttt{7}  & 7  & Route 7 Komatsugawa Line \\
\texttt{9}  & 9  & Route 9 Fukagawa Line \\
\texttt{10} & 10 & Route 10 Harumi Line \\
\texttt{11} & 11 & Route 11 Daiba Line \\
\bottomrule
\end{tabular}

\end{table}

\begin{table}[tbp]
  \caption{Summary statistics of the toy network data testbed and the line name codebook.}
  \label{tab:network_stats}
  \centering
  \small
  \begin{tabular}{lr@{\hspace{1.5em}}lr}
\toprule
\multicolumn{2}{l}{\textit{Topology}} & \multicolumn{2}{l}{\textit{Codebook}} \\
\midrule
Nodes                 & 623 & OD pairs (grid)            & 13,924 \\
Directed links        & 682 & Reachable OD pairs         & 10,094 \\
Named lines           & 14  & OD pairs with routes       & 10,094 \\
Junctions (JCT)       & 24  & Codewords $|\mathcal{C}|$  & 24,193 \\
Entry ramps (IN)      & 118 & Distinct $\Phi$            & 24,193 \\
Exit ramps (OUT)      & 118 & Vocabulary $|\mathcal{T}|$ & 274 \\
\midrule
\multicolumn{2}{l}{\textit{Length statistics}} & \multicolumn{2}{l}{\textit{Index size}} \\
\midrule
Median links / route             & 36   & Trie states            & 44,810 \\
Median tokens / route            & 15   & DAFSA states           & 3,269 \\
Median compression $|\Phi(r)|/|r|$ & 0.42 & DAFSA / Trie reduction & 92.7\% \\
\bottomrule
\end{tabular}

\end{table}

\subsection{Synthetic probe data}
\label{sec:setup_probe}

\begin{figure}[tbp]
  \centering
  \includegraphics[width=0.96\textwidth]{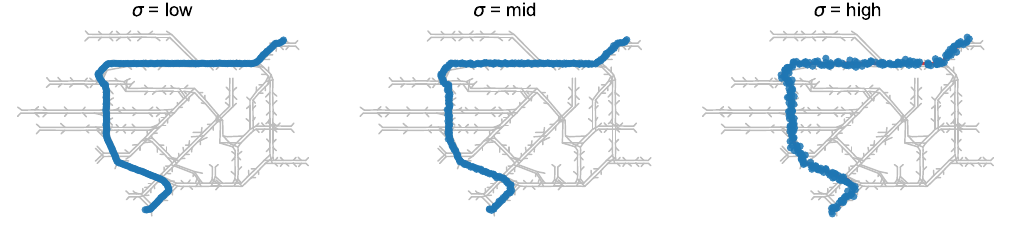}
  \caption{Probe spread on one long route at three relative noise levels (low / mid / high). Red: ground truth route polyline; blue dots: noisy probes.}
  \label{fig:probe_sigma}
\end{figure}

Ground truth routes are sampled uniformly from the codebook $\mathcal{C}$ (with the constraint $|r| \geq 5$ links to avoid degenerate cases), and synthetic probe trajectories are generated along the polyline implied by each route's node coordinate sequence.
The polyline is sampled at a fixed step relative to the typical link length, yielding raw probe points $\{p_{1}^{\circ}, \ldots, p_{T}^{\circ}\}$, and isotropic Gaussian noise is added at each point: $p_{t} = p_{t}^{\circ} + \varepsilon_{t}$ with $\varepsilon_{t} \sim \mathcal{N}(0, \sigma^{2} I)$.
Because the toy network coordinates are deformed, $\sigma$ is reported in \emph{link length units} (LLU): one LLU equals one typical link length, and we sweep $\sigma$ over a range that brackets realistic GPS noise on a single link.
\cref{fig:probe_sigma} shows how a single long route's probe trajectory scatters as $\sigma$ moves from low to high.

\begin{figure}[tbp]
  \centering
  \includegraphics[width=0.96\textwidth]{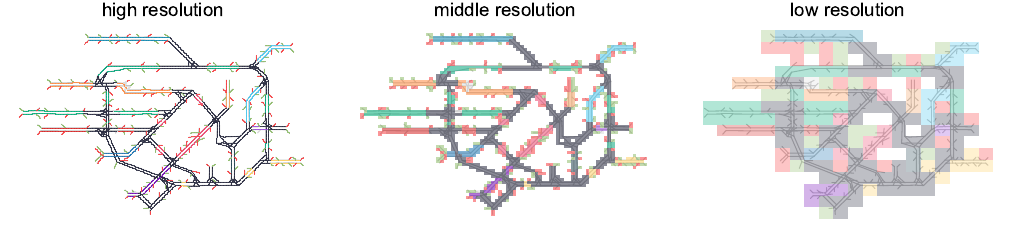}
  \caption{Mesh quantizer at three relative resolutions. Each cell is shaded by its dominant token; the network is overlaid in grey. Fine cells track individual line segments at the cost of empty cells under noise; coarse cells lump distinct lines together.}
  \label{fig:mesh_panels}
\end{figure}

The line assignment of each noisy probe is performed by the mesh quantizer; the per probe candidate token set is run length contracted on its argmax token to obtain the observation sequence fed to the decoder.
The mesh resolution is treated as a categorical factor with three settings (high, middle, low), whose boundary effects appear visually in \cref{fig:mesh_panels} and quantitatively in \cref{sec:res_mesh}.

\subsection{Evaluation metrics}
\label{sec:setup_metrics}

We evaluate decoding performance with three metrics at three granularities.
The \emph{exact match rate} (top 1) is the fraction of trajectories for which the decoded route is identical to the ground truth route as a link sequence; it captures global route level correctness.
The \emph{token F1} measures the agreement between the token sequences $\Phi(\hat{R})$ and $\Phi(R)$ of the decoded and ground truth routes (their multiset overlap, treating each token as one item); it captures token level correctness and is robust to the within OD ambiguity that the codebook representation deliberately equates.
The \emph{top 5 hit rate} is the fraction of trajectories for which the ground truth route lies within the top 5 returned codewords; it captures whether the correct route is in the soft decision shortlist.
Applications such as travel time estimation, origin destination inference, route based tolling, and bottleneck analysis require the exact match, so we adopt it as the primary metric, with token F1 and top 5 hit as diagnostics.
We do not report a separate link level F1 in this paper, since the codebook decoder returns a route in $\mathcal{R}$ by construction and the link level multiset is determined by the route choice; on a per query basis link F1 reduces to a thresholded variant of exact match.

\section{Results}
\label{sec:results}

\subsection{Experiment 1: accuracy under observation noise}
\label{sec:res_acc}

\begin{figure}[tbp]
  \centering
  \includegraphics[width=0.5\textwidth]{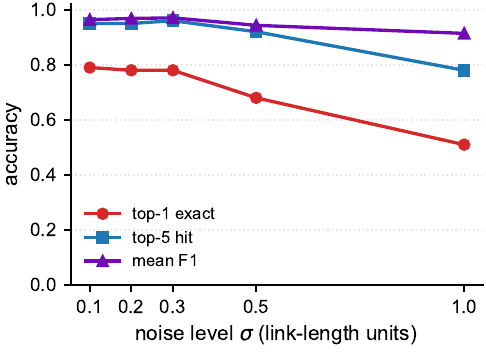}
  \caption{Experiment 1: top-1 exact, top-5 hit, and token $F_{1}$ as a function of noise $\sigma$ (link length units).}
  \label{fig:accuracy}
\end{figure}

The first experiment varies the observation noise level $\sigma$ over five values (in link length units, LLU; one LLU $= 280$\,m, the typical link length) at a fixed mesh resolution of 250\,m.
Each noise level is evaluated over 100 trials drawn from 50 distinct ground truth routes (with two probe noise seeds per route).
\cref{tab:main} reports the resulting top-1 exact, top-5 hit, and token level $F_{1}$ values; \cref{fig:accuracy} shows the same numbers as a function of $\sigma$.

\begin{table}[tbp]
  \caption{Experiment 1: decoding accuracy at five noise levels ($\sigma$ in link length units, LLU). Mesh resolution: middle.}
  \label{tab:main}
  \centering
  \small
  \begin{tabular}{rrrr}
\toprule
$\sigma$ (LLU) & top-1 & top-5 & mean F1 \\
\midrule
0.10 & 0.79 & 0.95 & 0.964 \\
0.20 & 0.78 & 0.95 & 0.968 \\
0.30 & 0.78 & 0.96 & 0.970 \\
0.50 & 0.68 & 0.92 & 0.943 \\
1.00 & 0.51 & 0.78 & 0.914 \\
\bottomrule
\end{tabular}
\end{table}

Three regimes appear in the accuracy curve (\cref{fig:accuracy}).
In the low to moderate noise regime ($\sigma \leq 0.30$\,LLU, i.e., up to about one third of the typical link length), the top-1 exact match rate is essentially flat at 0.78--0.79; the proposed method either decodes the correct route or returns a route in the same origin destination neighbourhood whose codeword is one substitution away.
In the same regime, the top-5 hit rate is at least 0.95: the ground truth route is almost always present in the top five candidates of the soft decision decoder.
In the moderate noise regime ($\sigma = 0.50$\,LLU), the top-1 rate falls to 0.68 while $F_{1}$ remains at 0.94; the per token alignment is still mostly correct, but the route level identity becomes harder to fix.
In the high noise regime ($\sigma = 1.00$\,LLU, equal to one full link length), the top-1 rate degrades to 0.51 and the top-5 rate to 0.78, but the token level $F_{1}$ stays above 0.91, which confirms that the line name representation is still informative about the route even when the exact codeword cannot be recovered.

The accuracy regime delineated above also identifies where the method does \emph{not} apply.
The framework targets sparse, structured expressway networks where every segment belongs to a named line and every junction has a name; on dense urban street networks without a named line structure no compact line name alphabet exists, and very short routes that collapse to one or two tokens lose codebook based disambiguation.
At very high noise the mesh quantizer's line assignment becomes close to uniform and the scoring loses discriminative power.
The accuracy gap between top-1 and top-5 originates in the within-OD ambiguity quantified in \cref{sec:res_coding}; when this within-OD separation matters, vocabulary enrichment (e.g., context tagged exit tokens such as \texttt{OUT@<IC>\_FROM\_<ring>}) widens it at the cost of index size.

\subsection{Experiment 2: mesh resolution sensitivity}
\label{sec:res_mesh}

\begin{figure}[tbp]
  \centering
  \includegraphics[width=0.96\textwidth]{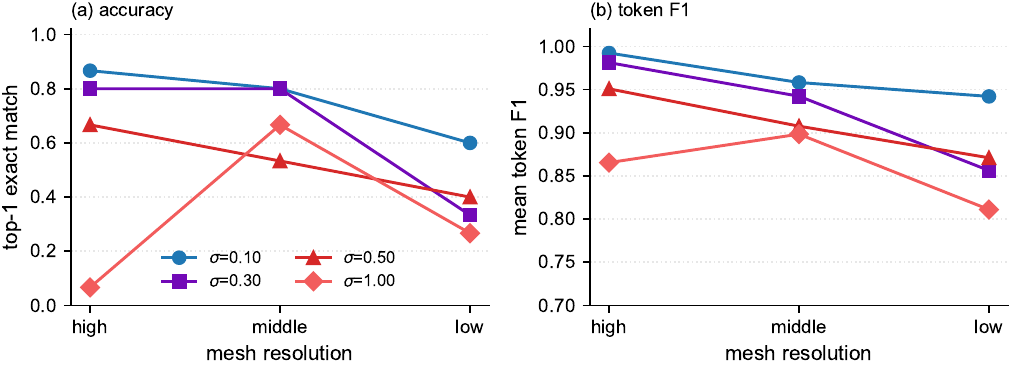}
  \caption{Experiment 2: (a) top-1 exact and (b) mean token $F_{1}$ at the three mesh resolution categories of \cref{fig:mesh_panels}, at four noise levels.}
  \label{fig:mesh_sensitivity}
\end{figure}

The mesh resolution controls the granularity of probe to line assignment.
Experiment 2 evaluates the three resolution categories (high, middle, low) of \cref{fig:mesh_panels} at four noise levels $\sigma \in \{0.10, 0.30, 0.50, 1.00\}$\,LLU.
At each cell of the $3 \times 4$ grid we evaluate 15 ground truth routes drawn from $\mathcal{C}$.
\cref{fig:mesh_sensitivity} shows the resulting top-1 accuracy and token level $F_{1}$ at each setting.

Three patterns appear in \cref{fig:mesh_sensitivity}.
First, at low noise ($\sigma = 0.10$\,LLU) the high and middle settings are within a few points of each other: every fine enough mesh assigns most probes to the correct cell.
Second, at high noise ($\sigma = 1.00$\,LLU) the high resolution setting collapses to a top-1 of 0.07: probes scatter into cells that the noisy displacement pushes them into, and the resulting observation sequence does not align with the true codeword.
The middle resolution setting holds top-1 above 0.4 at the same noise level because the cell size now absorbs the noise.
Third, at every noise level the low resolution setting underperforms: cells are large enough that the dominant token rule conflates lines that physically share the cell, and the observation sequence becomes ambiguous.
The token level $F_{1}$ in panel~(b) shows the same U shape with smaller amplitude: even when the route level identity is missed, the token level overlap remains above 0.85.
The middle setting balances the two failure modes and is the operating point used in the rest of the paper.

\subsection{Experiment 3: indexing structures and decoding speed}
\label{sec:res_index}

Experiment 3 compares four implementations of the same scored alignment decoder: brute force linear scan over the codebook, the trie variant, the DAFSA variant with path traced route lookup, and the DAFSA $\times$ Levenshtein product with a bounded edit budget.
All four share the cost model $(\mathrm{sub}, \mathrm{ins}, \mathrm{del}) = (1, 1, 1)$ and the symmetric Levenshtein distance over candidate sets; they differ only in the index structure on which the search runs.
The benchmark uses 20 ground truth routes at $\sigma = 0.30$\,LLU; each query is run on every variant and timed individually.

\begin{figure}[tbp]
  \centering
  \includegraphics[width=0.96\textwidth]{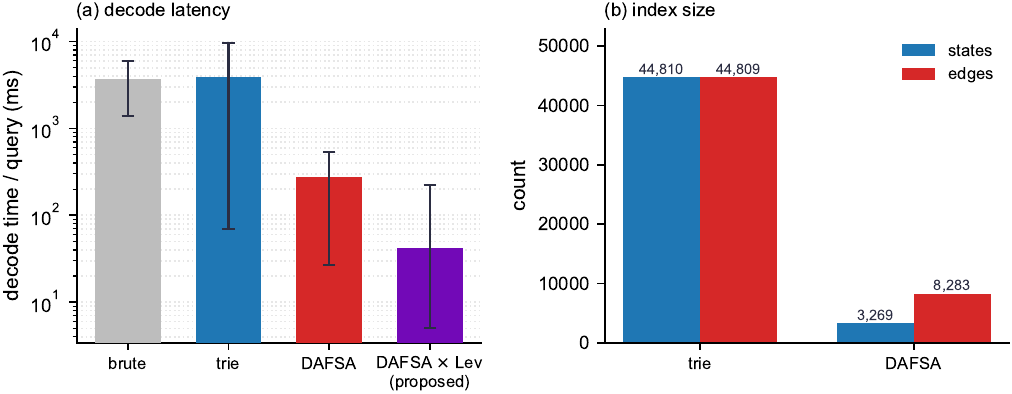}
  \caption{Experiment 3: (a) median per query decode time (whiskers: min--p95; log axis) and (b) index size for the four decoder variants. The DAFSA $\times$ Levenshtein bar (rightmost in panel (a)) is the proposed decoder.}
  \label{fig:index_compare}
\end{figure}

\cref{fig:index_compare} reports the median decode time per query for each variant alongside the index sizes; \cref{tab:exp03_index} tabulates the numerical breakdown.
\cref{tab:complexity} summarizes the corresponding theoretical complexity.

\begin{table}[tbp]
  \caption{Experiment 3: per query decode latency at $\sigma = 0.30$\,LLU, 20 queries. ``DAFSA $\times$ Lev'' is the proposed decoder.}
  \label{tab:exp03_index}
  \centering
  \small
  \begin{tabular}{lrrr}
\toprule
index / decoder & median (ms) & mean (ms) & 95\%ile (ms) \\
\midrule
brute & 3689.7 & 3956.0 & 6046.5 \\
trie & 3881.2 & 4302.5 & 9527.4 \\
DAFSA & 278.1 & 299.2 & 533.1 \\
DAFSA~$\times$~Lev & 41.6 & 64.1 & 225.1 \\
\midrule
trie states / edges & \multicolumn{3}{r}{ 44,810 / 44,809 } \\
DAFSA states / edges & \multicolumn{3}{r}{ 3,269 / 8,283 } \\
DAFSA / trie states  & \multicolumn{3}{r}{ 7.3\% (92.7\% reduction) } \\
\bottomrule
\end{tabular}
\end{table}

\begin{table}[tbp]
  \caption{Time complexity of the decoding variants. $|\mathcal{T}_{\mathcal{C}}|$: trie nodes; $|\mathcal{A}_{\mathcal{C}}|$: DAFSA states; $|\mathcal{A}_{\times}|$: product automaton size.}
  \label{tab:complexity}
  \centering
  \small
  \begin{tabular}{ll}
    \toprule
    Method & Time complexity \\
    \midrule
    Brute force DP & $O(|\mathcal{C}| \cdot L_{R} \cdot L_{\mathrm{obs}})$ \\
    Trie DP & $O(|\mathcal{T}_{\mathcal{C}}| \cdot L_{\mathrm{obs}})$ \\
    DAFSA DP & $O(|\mathcal{A}_{\mathcal{C}}| \cdot L_{\mathrm{obs}})$ \\
    DAFSA $\times$ Levenshtein (Dijkstra) \textbf{(proposed)} & $O(|\mathcal{A}_{\times}| \log |\mathcal{A}_{\times}|)$ \\
    \bottomrule
  \end{tabular}
\end{table}

The empirical ordering of decode latencies is the same as the theoretical ordering: brute force is the slowest, the trie is comparable to brute force at this codebook size, the DAFSA is roughly an order of magnitude faster than the trie, and the DAFSA $\times$ Levenshtein product is fastest by another order of magnitude.
The intuition is that suffix sharing in the DAFSA collapses the many same-OUT token suffixes of toy network data into a small number of accept states (24 distinct \texttt{JCT\_*} tokens, 118 distinct \texttt{OUT@*} tokens, against 24{,}193 codewords), so a single search visits far fewer states than a search on the trie.
The Levenshtein product additionally caps the search radius via the edit budget, which is the largest individual speed up at low to moderate noise.
None of the four runtimes depends directly on the number of links $N$: the cost is governed by the indexed codebook size and the observation length $L_{\mathrm{obs}}$.
The observation length itself stays short because consecutive repeated tokens are run length contracted (median codeword length is 15 tokens against a median of 36 links per route).

\subsection{Pairwise separation of the codebook}
\label{sec:res_coding}

The codebook on toy network data has $|\mathcal{C}| = 24{,}193$ codewords over an alphabet of $|\mathcal{T}| = 274$ tokens (14 \texttt{LINE\_*}, 24 \texttt{JCT\_*}, 118 \texttt{IN@*}, 118 \texttt{OUT@*}).
Within an origin destination pair, alternative routes typically differ in their interior junction sequence and ring traversal direction, with codeword length differences of a few tokens.
A representative pair of same-OD codewords (entry $\to$ exit at distinct interior junctions) is
\begin{align}
  &R_{a}: \texttt{[IN@Shibuya-Up, LINE\_3, JCT\_Tanimachi, LINE\_C1, OUT@Kasumigaseki-In]}, \notag \\
  &R_{b}: \texttt{[IN@Shibuya-Up, LINE\_3, JCT\_Ohashi, LINE\_C2, OUT@Kasumigaseki-In]},
\end{align}
which differ at the interior \texttt{JCT\_*} and \texttt{LINE\_*} tokens that distinguish the C1 versus C2 ring choice.
\par
Same origin to destination routes are the hardest region to separate, and inner loop versus outer loop traversals of the same ring are distinguished only through the interior junction sequence.
This is the rationale for moving the JCT identity from links to nodes (\cref{sec:setup_nets}): without that, two physically distinct C1 traversals collapse to the same codeword and become indistinguishable.
The scored alignment decoder of \cref{eq:decoder} is not a hard decision matcher; it accumulates continuous scores from the confusion matrix of \cref{eq:confusion} and the DP skip penalties, and the noise robustness reported in \cref{sec:res_acc} reflects this continuous scoring rather than any explicit error correction guarantee.

\section{Conclusion}
\label{sec:conclusion}

We presented a map matching method for urban expressway networks that operates at the line name granularity: routes become token sequences over a finite line and junction alphabet, the set of such sequences is a route codebook, probes are tokenised by a mesh quantizer (a precomputed coordinate to token grid), and decoding is scored alignment of the observation against codebook members.
By construction the decoder returns a member of the codebook, that is, a valid candidate route.
On a toy network data set, a deliberately deformed replica of the Tokyo Metropolitan Expressway, the method is robust at moderate GPS noise; suffix shared automaton indexing combined with a Levenshtein automaton yields a substantial decoding speed up over a linear scan; mesh resolution shows a U shaped sensitivity that bounds the operating range from both sides.
\par
The evaluation in this preprint is limited to synthetic probe trajectories on a topology with deliberately deformed coordinates, with hand set confusion matrix and skip penalty parameters and no head to head HMM comparison.
The primary remaining task before journal submission is therefore evaluation on real probe data collected on the Tokyo Metropolitan Expressway through established research channels such as ETC2.0; the same evaluation will also support data driven tuning of the channel model parameters and a structured comparison with an HMM baseline.

\section*{Acknowledgements}

The authors thank Hiroki Abe, Taisei Koshiji, Shuya Miyazu, and Seo Yanagawa for their help with organising the network and probe data used in this study.

\clearpage
\bibliographystyle{elsarticle-harv}
\bibliography{references}

\appendix

\section{Glossary of cross field terminology}
\label{app:glossary}

The paper uses vocabulary from several fields; the correspondences are listed for reference.
\cref{tab:glossary} summarizes the correspondences used throughout the text.

\begin{table}[tbp]
  \caption{Cross field terminology used in this paper.}
  \label{tab:glossary}
  \centering
  \small
  \begin{tabular}{@{}p{3.5cm} p{4.2cm} p{8.0cm}@{}}
    \toprule
    Field & Term & Meaning in this paper \\
    \midrule
    Codebook vocabulary
      & Codebook $\mathcal{C}$         & Set of token sequences of candidate routes, \cref{eq:codebook} \\
      & Codeword $\bm{c} = \Phi(R)$    & Token sequence of a single route \\
      & Vocabulary $\mathcal{T}$       & Finite set of line and junction tokens \\
      & Minimum distance $d_{\min}$    & Minimum Hamming distance in $\mathcal{C}$, \cref{eq:dmin} \\
    \midrule
    Formal language
      & Trie                           & Shared prefix tree index of $\mathcal{C}$ \\
      & DAFSA                          & Minimal acyclic DFA for $\mathcal{C}$ \\
      & Levenshtein automaton          & Accepter of strings within edit distance $k$ \\
      & Product automaton              & Decoding graph for scored alignment search \\
    \midrule
    Speech recognition analogue
      & WFST composition               & Structurally isomorphic to DAFSA $\times$ Levenshtein \\
    \midrule
    Transportation
      & Link $\ell$                    & Directed road segment \\
      & Line                           & Named roadway grouping consecutive links \\
      & Junction (JCT)                 & Transfer point between lines \\
      & Interchange (IC)               & Entry/exit to/from the expressway \\
      & OD pair                        & Origin destination pair \\
    \midrule
    Specific to this paper
      & Mesh quantizer                 & Probe to line lookup table, \cref{eq:quantizer} \\
      & Occurrence numbered token      & Line token with traversal index \texttt{Line.k} \\
      & Out of codebook route          & Link sequence returned by a per link decoder that is not in $\mathcal{R}$ \\
    \bottomrule
  \end{tabular}
\end{table}

\end{document}